\documentclass[11pt,a4paper]{article}
\usepackage[latin1]{inputenc}
\usepackage{amssymb}
\usepackage{pictex}
\usepackage{epsfig}
\newtheorem{rem}{Remark}[section]
\newtheorem{prop}{Proposition}[section]

\newtheorem{cor}{Corollary}[section]

\newcommand{\disp}{\displaystyle}
\newcommand{\bprf}{{\it Proof.~}}
\newcommand{\eprf}{\hfill $\square$ \bigskip\par}

\newcommand{\pitr}{ \mathbb{P}_3}

\newcommand{\oc}{\mathcal{O}}
\newcommand{\piu}{\mathbb{P}_1}

\newcommand{\ci}{ \mathbb{C}}
\newcommand{\er}{ \mathbb{R}}

\newcommand{\tu}{\tilde{T}\times 1}
\newcommand{\ii}{\tilde{I}\times 1}
\newcommand{\oo}{\tilde{O}\times 1}
\newcommand{\ttt}{1\times\tilde{T}}
\newcommand{\ooo}{1\times\tilde{O}}
\newcommand{\iii}{1\times\tilde{I}}

\begin{document}

\author{A.~Sarti
\\
FBR 17 Mathematik, Universit\"at Mainz\\
Staudingerweg 9, 55099 Mainz\\
Germany\\
e-mail:sarti@mathematik.uni-mainz.de\\}

\title{A geometrical construction for the polynomial invariants of
some reflection groups}

\maketitle

%\begin{center}
%{\Large A geometrical construction for the invariant polynomials  of some reflection groups}\\[0.3cm]
%{\rm \today}\\[0.3cm]
%A.Sarti\\[0.3cm]
%July 19, 2001\\[0.3cm]
%FBR 17 Mathematik, Universit\"at Mainz\\
%Staudingerweg 9, 55099 Mainz\\
%Germany\\
%e-mail:sarti@mathematik.uni-mainz.de\\
%\end{center}

\begin{abstract}
In these notes we investigate the ring of real polynomials in four variables, which are invariant under the action of the reflection groups $[3,4,3]$ and $[3,3,5]$.  It is well known that they are rationally generated in degree 2,6,8,12 and 2,12,20,30. We give a different proof of this fact by giving explicit equations for the generating polynomials.
\end{abstract}

\section{Introduction}
There are four groups generated by reflections 
which operate on the four-dimensional Euclidian space. These are the symmetry  groups of some regular four dimensional polytopes and are described in \cite{coxeter2} p. 145 and table I p. 292-295. With the notation there the groups and their orders are 

%\begin{center}

\renewcommand{\arraystretch}{1.3}

\begin{eqnarray*}
\begin{array}{c|cccc}
{\rm Group } & [3,3,3]&[3,3,4]&[3,4,3]&[3,3,5]\\
\hline
{\rm Order}  & 120&384&1152&14400\\
\end{array}
\end{eqnarray*}

%\end{center}

\renewcommand{\arraystretch}{1.0}

\noindent
They operate in a natural way on the ring of polynomials $R=\er[x_0,x_1,x_2,x_3]$ and it is well known that the ring of invariants $R^{G}$ ($G$  one of the groups above) is algebraically generated by a set of four independent polynomials (cf. \cite{burnside} p. 357).  Coxeter shows in \cite{coxeter1} that the rings $R^G$, $G= [3,3,3]$ or $[3,3,4]$ are generated in degree $2,3,4,5$ resp. $2,4,6,8$ and since the product of the degrees is equal to the order of the group, any other invariant polynomial is a combination with real coefficients of products of these invariants (i.e., in the terminology of \cite{coxeter1}, the ring $R^G$ is {\it rationally generated} by the polynomials). Coxeter also gives equations for the generators. In the case of the groups $[3,4,3]$ and $[3,3,5]$ he recalls a result of Racah, who shows with the help of the theory of Lie groups that the rings $R^G$ are rationally generated in degree $2,6,8,12$ resp. $2,12,20,30$ (cf. \cite{racah}). Neither Coxeter nor Racah give equations for the polynomials. In these notes we construct the generators and give a different proof of the result of Racah. 
%As remarked by Coxeter in \cite{coxeter1}, we have a trivial invariant in degree two, which is a positive definite quadratic form. 
The invariant of degree two is well known ( cf. \cite{coxeter1} ) and can be given as
\begin{eqnarray*}
q=x_0^2+x_1^2+x_2^2+x_3^2.
\end{eqnarray*}
We construct the other invariants in a completely geometrical way. For proving that our polynomials together with the quadric generate the ring $R^G$, we show some relations between them and the invariant forms of the binary tetrahedral group and of the binary icosahedral group.\\
%Moreover except for Propositions \ref{C1inv} and \ref{quad}, we do not use the explicit equations of the polynomials. 
%Observe that an equation for a $[3,4,3]$-invariant polynomial of degree six and for a $[3,3,5]$-invariant polynomial of degree twelve was already given by the author in \cite{sarti} by using the computer programm MAPLE.\\
It is a pleasure to thank W. Barth of the University of Erlangen for many helpful discussions.

\section{Notations and preliminaries}
Denote by $R$ the ring of polynomials in four variables with real coefficients $\er[x_0,x_1,x_2,x_3]$, by $G$ a finite group of homogeneous linear substitutions, and by $R^{G}$ the ring of invariant polynomials.\\
1. A set of polynomials $F_1,\ldots,F_n$  in $R$ is called {\it algebraically dependent} if there is a non trivial relation 
\begin{eqnarray*}
\sum \alpha_I (F_1^{i_1}\cdot\ldots \cdot F_n^{i_n})=0,
\end{eqnarray*}
where $I=(i_1,\ldots,i_n)\in\mathbb{N}^n$, $\alpha_I\in\er$.\\ 
%and $\sum d_1i_1+\ldots d_ni_n=m$.
2. The polynomials are called {\it algebraically independent} if they are not dependent. For the ring $R^G$, there always exists a set of four algebraically independent polynomials (cf. \cite{burnside}, thm. I, p. 357).\\
3. We say that $R^G$  is {\it algebraically generated} by a set of polynomials $F_1,\ldots,F_4$, if for any other polynomial $P\in R^G$ we have an algebraic relation 
\begin{eqnarray*}
\sum \alpha_I (P^{i_0} \cdot F_1^{i_1}\cdot\ldots\cdot F_4^{i_4})=0.
\end{eqnarray*}
4. We say that the ring $R^G$ is {\it rationally generated} by a set of polynomials $F_1,\ldots,F_4$, if for  any other polynomial $P\in R^G$ we have a relation 
\begin{eqnarray*}
\sum \alpha_I (F_1^{i_1}\cdot\ldots\cdot F_4^{i_4})=P,~~\alpha_I\in\mathbb{R}
\end{eqnarray*}
5. The four polynomials of 3 are called a {\it basic set} if they have the smallest possible degree (cf. \cite{coxeter1}).\\
%6. In \cite{coxeter1} Coxeter shows that if the four independent polynomials have the property that the product of their degrees is exactly $|G|$, then they are a basic set and they generate $R^G$ rationally.\\
6. There are two classical $2:1$ coverings 
\begin{eqnarray*}
\rho:SU(2)\rightarrow SO(3)~~~\hbox{and}~~~\sigma:SU(2)\times SU(2)\rightarrow SO(4),
\end{eqnarray*}
we denote by $T,O,I$ the tetrahedral group, the octahedral group and the icosahedral group in $SO(3)$ and by $\tilde{T}$, $\tilde{O}$, $\tilde{I}$ the corresponding binary groups in $SU(2)$ via the map $\rho$. The $\sigma$-images of $\tilde{T}\times\tilde{T}$, $\tilde{O}\times\tilde{O}$ and $\tilde{I}\times\tilde{I}$ in $SO(4)$ are denoted by $G_6$, $G_8$ and $G_{12}$. By abuse of notation we write $(p,q)$ for the image in $SO(4)$ of an element $(p,q)\in SU(2)\times SU(2)$. As showed in 
\cite{sarti} (3.1) p. 436, the groups $G_6$ and $G_{12}$ are subgroups of index four respectively two in the reflections groups $[3,4,3]$ and $[3,3,5]$.

\section{Geometrical construction}
%In \cite{sarti} we have seen that $[3,4,3]$ contains $G_6$ the be-polyhedral tetrahedral group and $[3,3,5]$ contains $G_{12}$ the be-polyhedral icosahedral group, which are subgroups of index $4$, resp. $2$. The matrices are elements $(p,q)\in SO(4)$ (richiamare). 
Denote by $\tilde{G}$ one of the groups $\tilde{T}$, $\tilde{O}$ or $\tilde{I}$. Clearly, the subgroups  $\tilde{G}\times 1$ and $1\times\tilde{G}$ of $SO(4)$ are isomorphic to $\tilde{G}$. Moreover, each of them operates on one of the two rulings of the quadric $\piu\times\piu$ and leaves invariant the other ruling (as shown in \cite{sarti}). We recall the lengths of the orbits of points under the action of the groups $T$, $O$ and $I$

%\begin{center}

\renewcommand{\arraystretch}{1.3}

\begin{eqnarray*}
\begin{array}{c|c|c|c}
{\rm group}&{\rm T } &{\rm O } &{\rm I }\\
\hline
{\rm lengths~of~the~orbits}&12,~6,~4&24,~12,~8,~6&60,~30,~20,~12\\
\end{array}
\end{eqnarray*}

%\end{center}

\renewcommand{\arraystretch}{1.0}

\noindent
These lines are fixed by elements $(p,1)\in\tilde{G}\times 1$ on one ruling, resp. $(1,p')\in1\times\tilde{G}$ on the other ruling of the quadric. Recall that these elements have two lines of fix points with eigenvalues $\alpha$, $\bar{\alpha}$ which are in fact the eigenvalues of $p$ and $p'$. We call two lines $L,L'$ of $\piu\times\piu$ a {\it couple} if $L$ is fixed by $(p,1)$ with eigenvalue $\alpha$ and $L'$ is fixed by $(1,p)$ with the same eigenvalue. 
\subsection{The invariant polynomials of $G_6$ and of $G_{12}$}
Consider the six couples of lines $L_1,L_1', \ldots, L_6, L_6'$ in $\piu\times\piu$ which form one orbit under the action of $\tu$, resp. $\ttt$, and denote by $f_{11}^{(6)},\ldots,f_{66}^{(6)}$ the six planes generated by such a couple of lines (and by abuse of notation their equation, too). Now set
\begin{eqnarray*}
F_6=\sum_{g\in\tu} g(f_{11}^{(6)}\cdot f_{22}^{(6)}\cdot\ldots\cdot f_{66}^{(6)})=\sum_{g\in\tu}g(f_{11}^{(6)})\cdot g(f_{22}^{(6)})\cdot\ldots\cdot g(f_{66}^{(6)}).
\end{eqnarray*}
Observe that an element $g\in\tu$ leaves each line of one ruling invariant and operates on the six lines of the other ruling. A similar action is given by an element of $\ttt$. Since we sum over all the elements of $\tu$, the action of  $1\times\tilde{T}$ does not give anything new, hence $F_6$ is $G_6$-invariant. Furthermore observe that $F_6$ has real coefficients. In fact, in the above product, for each plane generated by the lines $L_i$, $L_i'$ we also take the plane generated by the lines which consist of the conjugate points. The latter has equation $\bar{f_{ii}}^{(6)}$, i.e., we have an index $j\not= i$ with $f_{jj}^{(6)} =\bar{f_{ii}}^{(6)}$  and the products $f_{ii}^{(6)}\cdot\bar{f_{ii}}^{(6)}$ have real coefficients.\\
Consider now the orbits of lengths eight and twelve under the action of $\oo$ and $\ooo$ and the planes $f_{ii}^{(8)}$, $f_{jj}^{(12)}$ generated by the eight, respectively by the twelve couples of lines. As before the polynomials
\begin{eqnarray*}
\begin{array}{lll}
F_8&=&\displaystyle\sum_{g\in\tu} g(f_{11}^{(8)}\cdot\ldots\cdot f_{88}^{(8)}),\\
&&\\
F_{12}&=&\displaystyle \sum_{g\in\tu} g(f_{11}^{(12)}\cdot\ldots\cdot f_{1212}^{(12)})
\end{array}
\end{eqnarray*}
%Since we sum over all the elements of $\tu$ , $F_8$ and $F_{12}$ 
are $G_6$-invariant and have real coefficients.\\
%\subsection{The invariants of $\gd$}
Finally we consider the lines of $\piu\times\piu$ which form orbits of length $12,20$ and $30$ under the action of $\ii$ resp. $\iii$. The planes generated by the couples of lines produce the $G_{12}$-invariant real polynomials  
\begin{eqnarray*}
\begin{array}{lll}
\Gamma_{12}&=&\displaystyle \sum_{g\in\ii} g(h_{11}^{(12)}\cdot\ldots\cdot h_{1212}^{(12)}),\\
&&\\
\Gamma_{20}&=&\displaystyle \sum_{g\in\ii} g(h_{11}^{(20)}\cdot\ldots\cdot h_{2020}^{(20)}),\\
&&\\
\Gamma_{30}&=&\displaystyle \sum_{g\in\ii} g(h_{11}^{(30)}\cdot\ldots\cdot h_{3030}^{(30)}).\\
\end{array}
\end{eqnarray*}
%are $G_{12}$-invariant. As in the case of $F_6$ all these polynomials have real coefficients.

\subsection{The invariant polynomials of the reflection groups}
We consider the matrices
\begin{eqnarray*}
C=\left( \begin{array} {cccc} 
1& 0&0&0\\
0& -1&0&0\\
0&0&-1&0\\
0&0&0&-1
\end{array} \right),~~
C'=\left( \begin{array} {cccc} 
1& 0&0&0\\
0& 1&0&0\\
0&0&0&1\\
0&0&1&0
\end{array} \right),
\end{eqnarray*}
 as in \cite{sarti} (3.1) p. 436, the groups generated by $G_6$, $C$, $C'$ and $G_{12}$, $C$ are the reflections groups $[3,4,3]$ respectively $[3,3,5]$.
\begin{prop}\label{C1inv}
1. The polynomials $F_6$, $F_8$, $F_{12}$, $\Gamma_{12}$, $\Gamma_{20}$, $\Gamma_{30}$ are $C$ invariant.\\
2. The polynomials $F_6$, $F_8$, $F_{12}$ are  $C'$ invariant.
\end{prop}
\bprf
1. The matrix $C$ interchanges the two rulings of the quadric, hence the polynomials $F_i$ and $\Gamma_j$ are invariant by construction. We prove 2 by a direct computation in the last section.\eprf
From this fact we obtain
\begin{cor}
The polynomials $q,F_6, F_8, F_{12}$ are $[3,4,3]$-invariant and the polynomials $q,\Gamma_{12},\Gamma_{20},\Gamma_{30}$ are $[3,3,5]$-invariant.
\end{cor}
Here we denote by $q$ the quadric $\piu\times\piu$.

%In section 4 by a direct computation we prove also
%\begin{prop}\label{quad}
%The quadric $q$ does not divide the polynomials $F_i$, $\Gamma_j$. Moreover $F_6$ does not divide $F_{12}$.
%\end{prop}
%\bprf It is an explicit computation which we give in the section \ref{comput}.
%Consider the points $p_1=(i\sqrt{2}:1:1:0)$ and $p_2=(1:i:0:0)$ of $\pitr$ then $q(p)=0$. Take now for each polynomial just one of the products of the planes, i.e. those given in the last section.
%Computing the value of this product on $g\cdot p_i$, $g\in(\tilde{G},1)$ and doing the sum we see that $F_6(p_1)\not=0$ $F_8(p_1)\not=0$ and $F_{12}(p_2)\not=0$. For the exact computation see the last section. In a similar way one shows the assertion in the remaining cases. 
%$\Gamma_j(p)\not=0$. A similar calculation shows that for $p=...$ we have $F_6(p)=0$, but $F_{12}(p)\not=0$.
%\eprf 
%All the polynomials $F_6,F_8,F_{12}, \Gamma_{12},\Gamma_{20}, \Gamma_{30}$ are invariant under the matrix $C=....$ because it interchanges just the lines of the two rulings (forse e' meglio mettere altri invarianti...). Hence $\Gamma_{12},\Gamma_{20}, \Gamma_{30}$ are $[3,3,5]$-invariant too (recall that $[3,3,5]=<G_{12},C>$). A direct computation (cf. section...) shows that $F_6,F_8,F_{12}$ are $C'=..$ invariant too.
\section{The rings of invariant forms}
%\section{The map $\mathcal{O}_{\mathbb{P}^3}(n)\longrightarrow \mathcal{O}_q(n,n)$}
%{\it Claim 0. The quadric $q$ does not divide $F_6,F_8,F_{12}, \Gamma_{12},\Gamma_{20}, \Gamma_{30}$ and $F_6$ does not divide $F_{12}$}\\
%\bprf
%It is a direct computation (cf. section...) because it is possible to find a point of $q$ which is not on $F_6,F_8,F_{12}, \Gamma_{12},\Gamma_{20}, \Gamma_{30}$ and a point of $F_6$ not in $F_{12}$ (migliorare).\eprf
Identify $\pitr$ with $\mathbb{P}M(2\times 2,\mathbb{C})$ by the map 
\begin{eqnarray}\label{identif}
\begin{array}{lll}
(x_0:x_1:x_2:x_3)&\mapsto&\left(\begin{array}{cc}
x_0+ix_1&x_2+ix_3\\
-x_2+ix_3&x_0-ix_1
\end{array}\right).
\end{array}
\end{eqnarray}
Furthermore consider the map 

\begin{eqnarray}\label{tensor}
\begin{array}{lll}
 \ci^2 \times \ci^2 &\longrightarrow & M(2\times 2,\mathbb{C})\\
((z_0,z_1),(z_2,z_3)) &\longmapsto & %\left(\begin{array}{c}
%z_0\\
%z_1
%\end{array}\right)\cdot \left(\begin{array}{cc}
%z_2& z_3
%\end{array}\right)=
\left( \begin{array} {cc} 
z_0z_2& z_0z_3\\
z_1z_2& z_1z_3
\end{array} \right)=\mathcal{Z}.
\end{array}
\end{eqnarray}
Then $\mathcal{Z}$ is a matrix of determinant $x_0^2+x_1^2+x_2^2+x_3^2=0$ which is the equation of $q$. Now denote by $\oc_{\pitr}(n)$  the sheaf of regular functions of degree $n$ on $\pitr$ and by $\oc_{q}(n,n)$ the sheaf of regular function of be-degree $(n,n)$ on the quadric $q$. We obtain a surjective map between the global sections 
\begin{eqnarray}\label{proj}
\begin{array}{llll}
\phi:&H^0(\oc_{\pitr}(n))&\longrightarrow &H^0(\oc_{q}(n,n))\\
\end{array}
\end{eqnarray}
by doing the substitution 
\begin{eqnarray*}
\begin{array}{ll}
x_0=\frac{\disp z_0z_2+z_1z_3}{\disp 2},&x_1=\frac{\disp z_0z_2-z_1z_3}{\disp 2i},\\
x_2=\frac{\disp z_0z_3-z_1z_2}{\disp 2}, &x_3=\frac{\disp z_0z_3+z_1z_2}{\disp 2i}
\end{array}
\end{eqnarray*}
in a polynomial $p(x_0,x_1,x_2,x_3)\in H^0(\oc_{\pitr}(n))$. Observe that $\phi(q)=0$.
%and identify $\pitr$ with $\mathbb{P}M(2\times 2,\mathbb{C})$ by the map
%\begin{eqnarray}\label{identif}
%\begin{array}{lll}
%(x_0:x_1:x_2:x_3)&\mapsto&\left(\begin{array}{cc}
%x_0+ix_1&x_2+ix_3\\
%-x_2+ix_3&x_0-ix_1
%\end{array}\right).
%\end{array}
%\end{eqnarray}
%In this way we get a surjective map
%\begin{eqnarray*}
%\begin{array}{lll}
%\pitr&\longrightarrow&\piu\times\piu\\
%(x_0:x_1:x_2:x_3)&\mapsto &\displaystyle (\frac{z_0z_2+z_1z_3}{2}:\frac{z_0z_2-z_1z_3}{2i}:\frac{z_0z_3-z_1z_2}{2}:\frac{z_0z_3+z_1z_2}{2i})\\
%\end{array}
%\end{eqnarray*}
%which induces a surjective map 
%\begin{eqnarray}\label{proj}
%\begin{array}{llll}
%\phi:&\oc_{\pitr}(n)&\longrightarrow &\oc_{q}(n,n)\\
%\end{array}
%\end{eqnarray}
%between the sheaf of regular functions of degree $n$ on $\pitr$ and the sheaf of regular function of be-degree $(n,n)$ on the quadric $q$. Observe that $\phi(q)=0$. \\
Now let
\begin{eqnarray*}
\begin{array}{lll}
t&=&z_0z_1(z_0^4-z_1^4),\\
W&=&z_0^8+14z_0^4z_1^4+z_1^8,\\
\chi&=&z_0^{12}-33(z_0^8z_1^4+z_0^4z_1^8)+z_1^8
\end{array}
\end{eqnarray*}
denote the $\tilde{T}$-invariant polynomials of degree $6,8$ and $12$  and 
let
\begin{eqnarray*}
\begin{array}{lll}
f&=&z_0z_1(z_0^{10}+11z_0^5z_1^5-z_1^{10}),\\
H&=&-(z_0^{20}+z_1^{20})+228(z_0^{15}z_1^5-z_0^5z_1^{15})-494z_0^{10}z_1^{10},\\
\mathcal{T}&=&(z_0^{30}+z_1^{30})+522(z_0^{25}z_1^{5}-z_0^5z_1^{25})-10005(z_0^{20}z_1^{10}+z_0^{10}z_1^{20})\\
\end{array}
\end{eqnarray*}
be the $\tilde{I}$-invariant polynomials of degree $12,20,30$ given by Klein in \cite{klein} p. 51-58. Put $t_1=t(z_0,z_1),~t_2=t(z_2,z_3)$, $W_1=W(z_0,z_1),~W_2=W(z_2,z_3)$ and analogously for the other invariants.\\
\begin{prop}\label{invpol}
If $p\in H^0(\mathcal{O}_{\mathbb{P}^3}(n))$ is $G_6$-invariant,  then:\\
\begin{eqnarray*}
\phi(p)=\sum _{I} \alpha_{I} t_1^{\alpha_1}t_2^{\alpha_1'}W_1^{\alpha_2}W_2^{\alpha_2'}\chi_1^{\alpha_3}\chi_2^{\alpha_3'}
\end{eqnarray*} 
If $p$ is $G_{12}$-invariant, then:
\begin{eqnarray*}
\phi(p)=\sum _{J} \beta_{J} f_1^{\beta_1}f_2^{\beta_1'}H_1^{\beta_2}H_2^{\beta_2'}\mathcal{T}_1^{\beta_3}\mathcal{T}_2^{\beta_3'}
\end{eqnarray*} 
where 
\begin{eqnarray*}
I=\{(\alpha_1,\alpha_1',\alpha_2,\alpha_2',\alpha_3, \alpha_3')|\alpha_i,\alpha_i'\in \mathbb{N}, 6\alpha_1+8\alpha_2+12\alpha_3=n,~6\alpha_1'+8\alpha_2'+12\alpha_3'=n \},\\
J=\{(\beta_1,\beta _1',\beta_2,\beta_2',\beta_3, \beta_3')|\beta_i,\beta_i'\in \mathbb{N}, 12\beta_1+20\beta_2+30\beta_3=n,~12\beta_1'+20\beta_2'+30\beta_3'=n \}.\\
\end{eqnarray*}
\end{prop}
%2) $f_1,f_2$ are $\tilde(G)$-invariant\\
%3) $f_1(x,y)$, $f_2(x,y)$ are equal up to scalar factor\\}
\bprf 
Put
\begin{eqnarray*}
\phi(p)=p'(z_0,z_1,z_2,z_3).
\end{eqnarray*}
An element  $g=(g_1,g_2)$ in $G_6$ or $G_{12}$ operates on $(x_0:x_1:x_2:x_3)\in\pitr$ by the matrix multiplication
\begin{eqnarray*}
g_1\left(\begin{array}{cc}
x_0+ix_1&x_2+ix_3\\
-x_2+ix_3&x_0-ix_1
\end{array}
\right)g_2^{-1}
\end{eqnarray*}
and on the matrix $\mathcal{Z}$ of (\ref{tensor}) by
\begin{eqnarray*}
\begin{array}{lll}
g_1\left(\begin{array}{cc}
z_0z_2&z_0z_3\\
z_1z_2&z_1z_3
\end{array}
\right)g_2^{-1}&=&g_1\left(\begin{array}{c}
z_0\\
z_1
\end{array}\right)\cdot \left(\begin{array}{cc}
z_2& z_3
\end{array}\right)g_2^{-1}.
\end{array}
\end{eqnarray*}
Clearly if $p$ is $G_6$- or $G_{12}$-invariant then also the projection $\phi(p)$ with the previous operation is. In particular for $g=(g_1,1)$ in $\tu$, resp. in $\ii$ the polynomial $p'$ is $\tu$-, respectively $\ii$-invariant as polynomial in the coordinates $(z_0:z_1)\in\piu$ and for any $(z_2:z_3)\in\piu$. On the other hand for $g=(1,g_2)$ in $\ttt$, resp. in $\iii$  the polynomial $p'$ is $\ttt$-, respectively $\iii$-invariant as polynomial in the coordinate $(z_2:z_3)\in\piu$ and for any $(z_0:z_1)\in\piu$. Hence $p'$ must be in the form of the statement.

\eprf
By a direct computation in section 4 we prove the following 
\begin{prop}\label{quad}
The quadric $q$ does not divide the polynomials $F_i$, $\Gamma_j$. Moreover, $F_6$ does not divide $F_{12}$.
\end{prop}
% Whenevr $\tilde{G}=\aq,~\ac$ we know the invariant polynomials under the action of the binary groups (cf. klein). In particular the only possibitlity is that $f''$ ia a product of two $\tilde{G}$-invariant polynomials i.e. $f''=f''_1(z_0,z_1)\cdot f''_2(z_2,z_3)$. Thos shows 1) and 2). it remains to prove 3). \eprf
\begin{cor}\label{phi}
We have $\phi(q)=0$, $\phi(F_6)=t_1\cdot t_2$, $\phi(F_8)=W_1\cdot W_2$, $\phi(F_{12})=\chi_1\cdot\chi_2$, $\phi(\Gamma_{12})=f_1\cdot f_2$, $\phi(\Gamma_{20})=H_1\cdot H_2$, $\phi(\Gamma_{30})=T_1\cdot T_2$ (up to some scalar factor).
\end{cor}
\bprf
This follows from Proposition \ref{invpol} and \ref{quad}
%{\it Claim 0} and {\it Claim 1}.\eprf
\begin{prop}\label{algebr}
The polynomials $q, F_6,F_8,F_{12}$, resp.  $q, \Gamma_{12},\Gamma_{20}, \Gamma_{30}$ are algebraically independent.\end{prop}
\bprf
Let $\sum_{I}\alpha_{I} q^{i_1} F_6^{i_2}F_8^{i_3}F_{12}^{i_4}=0$ and $\sum_{J}\beta_{J} q^{j_1} \Gamma_{12}^{j_2}\Gamma_{20}^{j_3}\Gamma_{30}^{j_4}=0$ be algebraic relations,  $I=(i_1,i_2,i_3,i_4)\in\mathbb{N}^4$, $J=(j_1,j_2,j_3,j_4)\in\mathbb{N}^4$, $\alpha_{I}, \beta_{J}\in\er$, then
\begin{eqnarray}\label{algrel}
\begin{array}{lll}
0&=&\phi(\sum_{I}\alpha_{I} q^{i_1} F_6^{i_2}F_8^{i_3}F_{12}^{i_4})\\
&=&\sum_{I'}\alpha_{I'} \phi(F_6)^{i_2'}\phi(F_8)^{i_3'}\phi(F_{12})^{i_4'}\\
&=&\sum_{I'}\alpha_{I'} t_1^{i_2'}t_2^{i_2'}W_1^{i_3'}W_2^{i_3'}\chi_1^{i_4'}\chi_2^{i_4'}
\end{array}
\end{eqnarray}
similarly
\begin{eqnarray}\label{algrel1}
\begin{array}{lll}
0&=&\phi(\sum_{J}\beta_{J} q^{j_1} \Gamma_{12}^{j_2}\Gamma_{20}^{j_3}\Gamma_{30}^{j_4})\\
&=&\sum_{J'}\beta_{J'} \phi(\Gamma_{12})^{j_2'}\phi(\Gamma_{20})^{j_3'}\phi(\Gamma_{30})^{j_4'}\\
&=&\sum_{J'}\beta_{J'} f_1^{j_2'}f_2^{j_2'}H_1^{j_3'}H_2^{j_3'}\mathcal{T}_1^{j_4'}\mathcal{T}_2^{j_4'}.
\end{array}
\end{eqnarray}
%Here $J=(j_2,j_3,j_4)\in\mathbb{N}^3$.
If the polynomials $t_1,W_1,\chi_1$ are fixed, we obtain a relation between $t_2,W_2$ and $\chi_2$, which is the same relation as for $t_1$, $W_1$ and $\chi_1$ if we fix $t_2,W_2$ and $\chi_2$. The same holds for the polynomials $f_1, H_1, \mathcal{T}_1$ and $f_2, H_2, \mathcal{T}_2$. From \cite{klein} p. 55  and p. 57 there are only the relations
\begin{eqnarray*}
108\,t_1^4-W_1^3+\chi_1^2=0,~~~~108\,t_2^4-W_2^3+\chi_2^2=0
\end{eqnarray*}
and 
\begin{eqnarray*}
\mathcal{T}_1^2+H_1^3-1728 f_1^5=0,~~~~\mathcal{T}_2^2+H_2^3-1728 f_2^5=0
\end{eqnarray*}
 between these polynomials. By multiplying these relations, however, it is not possible to obtain expressions like (\ref{algrel}) and (\ref{algrel1}).
\eprf
\begin{cor}
 The polynomials $q, F_6,F_8,F_{12}$, resp.  $q, \Gamma_{12},\Gamma_{20}, \Gamma_{30}$ generate rationally the ring of invariant polynomials of $[3,4,3]$, resp. $[3,3,5]$. \end{cor}
\bprf (cf. \cite{coxeter1} p. 775)
By Proposition \ref{algebr} and Proposition \ref{quad} these are algebraically independent, moreover the products of their degrees are
\begin{eqnarray*}
2\cdot 6\cdot 8\cdot 12=1152~~\mbox{and}~~2\cdot 12\cdot20\cdot 30=14400,
\end{eqnarray*}
which are equal to the order of the groups $[3,4,3]$ and $[3,3,5]$. By \cite{coxeter1} this implies the assertion.\eprf
\section{Explicit computations}\label{comput}
We recall the following matrices of $SO(4)$ (cf. \cite{sarti}).

\begin{eqnarray*}
\begin{array}{rr}
%(q_1,1)=
%\left( \begin{array} {cccc}
%0 & -1& 0 & 0 \\
%1 & 0& 0 & 0 \\
%0 & 0& 0 & -1 \\
%0 & 0& 1 & 0 
%\end{array}\right),&

(q_{2},1)=
\left( \begin{array} {cccc}
0 & 0& -1 & 0 \\
0 & 0& 0 & 1 \\
1 & 0& 0 & 0 \\
0 & -1& 0 & 0 
\end{array} \right),&
(1,q_{2})=
\left( \begin{array} {cccc}
0 & 0& 1 & 0 \\
0 & 0& 0 & 1 \\
-1 & 0& 0 & 0 \\
0 & -1& 0 & 0 
\end{array} \right),
\end{array}
\end{eqnarray*}
%\begin{eqnarray*}
%\begin{array}{ll}
%(1,q_{1})=
%\left( \begin{array} {cccc}
%0 & 1& 0 & 0 \\
%-1 & 0& 0 & 0 \\
%0 & 0& 0 & -1 \\
%0 & 0& 1 & 0 
%\end{array} \right),&
%\end{array}
%\end{eqnarray*}
\begin{eqnarray*}
\begin{array}{ll}
(p_3,1)=\frac{1}{2}
\left( \begin{array} {cccc}
1& -1& 1& -1 \\
1 & 1& -1 & -1 \\
-1 & 1& 1 & -1 \\
1 & 1& 1 & 1 
\end{array} \right),&
(1,p_3)=\frac{1}{2}
\left( \begin{array} {cccc}
1 & 1& -1 & 1 \\
-1 & 1& -1 & -1 \\
1 & 1& 1 & -1 \\
-1 & 1& 1 & 1 
\end{array} \right),\\
\end{array}
\end{eqnarray*}

\begin{eqnarray*}
\begin{array}{ll}
(p_4,1)=\frac{1}{\sqrt{2}}
\left( \begin{array} {cccc}
1& -1& 0& 0 \\
1 & 1& 0 & 0 \\
0 & 0& 1 & -1 \\
0 & 0& 1 & 1 
\end{array} \right),&
(1, p_4)=\frac{1}{\sqrt{2}}
\left( \begin{array} {cccc}
1 & 1& 0 & 0 \\
-1 & 1& 0& 0 \\
0 & 0& 1 & -1 \\
0 & 0& 1 & 1 
\end{array} \right),
\end{array}
\end{eqnarray*}
\begin{eqnarray*}
\begin{array}{c}
(p_5,1)=\frac{1}{2}
\left( \begin{array} {cccc}
\tau & 0& 1-\tau & -1 \\
0 & \tau & -1 & \tau-1 \\
\tau-1 & 1& \tau & 0 \\
1 & 1-\tau& 0 & \tau 
\end{array} \right),\\
(1,p_5)=\frac{1}{2}
\left( \begin{array} {cccc}
\tau & 0& \tau-1 & 1 \\
0 & \tau& -1 & \tau-1 \\
1-\tau & 1& \tau & 0 \\
-1 & 1-\tau & 0 & \tau 
\end{array} \right),
\end{array}
\end{eqnarray*}
%$\sigma_i^2=\pi_j^j=\pi_j'^j=-1,~i=1,2,3,4;~j=3,4,5.$, and 
where $\tau=\frac{1}{2}(1+\sqrt{5})$. Then we have

\begin{center}

\renewcommand{\arraystretch}{1.3}

\begin{eqnarray*}
\begin{array}{c|c}
{\rm Group } &{\rm Generators } \\
\hline
G_6& (q_2,1), (1,q_2), (p_3,1), (1,p_3)\\
G_8& (q_2,1), (1,q_2), (p_3,1), (1,p_3), (p_4,1), (1,p_4)\\
G_{12}& (q_2,1), (1,q_2), (p_3,1), (1,p_3), (p_5,1), (1,p_5)
\end{array}
\end{eqnarray*}

\end{center}

%\renewcommand{\arraystretch}{1.0} 

 %The $G_6$ is generated by $\sigma_1,\sigma_2,\sigma_3,\sigma_4,\pi_3,\pi_3'$, $G_8$ by $\sigma_1,\sigma_2,\sigma_3,\sigma_4,\pi_3,\pi_3',\pi_4,\pi_4'$ and $G_{12}$ by $\sigma_1,\sigma_2,\sigma_3,\sigma_4,\pi_5,\pi_5'$.\\
Now we can write down the equations of the fix lines on $\piu\times\piu$ and those of the planes which are generated by a couple of lines. The products of planes of section 2.1 in the case of the group $G_6$ are 
%In this section we write down the expression of the product of planes, which we consider in section...... In fact it  is possible to choose these planes in different ways. In any case one gets the same result 
\begin{eqnarray*}
\begin{array}{lll}
f_{11}^{(6)}\cdot f_{22}^{(6)}\cdot\ldots\cdot f_{66}^{(6)}&=&(x_2-ix_3)(x_1+ix_3)(x_2+ix_3)(x_1-ix_2)(x_1-ix_3)(x_1+ix_2),\\
f_{11}^{(8)}\cdot f_{22}^{(8)}\cdot\ldots\cdot f_{88}^{(8)}&=&(x_1+ax_2-bx_3)(x_1+bx_2-ax_3)(x_1-ax_2-bx_3)(x_1-ax_3-bx_2)\\
&&(x_2+bx_1-ax_3)(x_2+ax_1-bx_3)(x_2-bx_1+ax_3)(x_2+bx_3-ax_1),\\
f_{11}^{(12)}\cdot f_{22}^{(12)}\cdot\ldots\cdot f_{1212}^{(12)}&=&(x_3-x_1+cx_2)(x_3-x_1-cx_2)(x_2+x_3-cx_1)(x_2+x_3+cx_1)\\
&&(x_3-x_2+cx_1)(x_3-x_2-cx_1)(x_1+x_2+cx_3)(x_1+x_2-cx_3)\\
&&(x_1+x_3-cx_2)(x_1+x_3+cx_2)(x_1-x_2+cx_3)(x_1-x_2-cx_3),
\end{array}
\end{eqnarray*}
with $a=(1/2)(1+i\sqrt{3}), b=(1/2)(1-i\sqrt{3}), c=i\sqrt{2}$.\\
%h_{11}^{(12)}\cdot\ldots\cdot h_{1212}^{(12)}\\
%h_{11}^{(20)}\cdot\ldots\cdot h_{2020}^{(20)})\\
%h_{11}^{(30)}\cdot\ldots\cdot h_{3030}^{(30)})\\
Then the $G_6$-invariant polynomials $F_6$, $F_8$ and $F_{12}$ have the following expressions
\begin{eqnarray*}
\begin{array}{lll}
F_6&=&x_0^6+x_1^6+x_2^6+x_3^6+5x_0^2x_1^2(x_0^2+x_1^2)+5x_1^2x_3^2(x_1^2+x_3^2)+5x_1^2x_2^2(x_1^2+x_2^2)\\
&&+6x_0^2x_2^2(x_0^2+x_2^2)+6x_0^2x_3^2(x_0^2+x_3^2)+6x_3^2x_2^2(x_2^2+x_3^2)+2x_0^2x_2^2x_3^2,\\
&&\\
F_8&=&3\sum x_i^8+12\sum x_i^6x_j^2+30\sum x_i^4x_j^4+24\sum x_i^4x_j^2x_k^2+144x_0^2x_1^2x_2^2x_3^2,\\
&&\\
%24x0^2x2^4x1^2+24x0^4x2^2x1^2+24x1^4x2^2x3^2+24x1^4x2^2x0^2\\+
%24x0^2x3^4x1^2+24x0^2x3^2x1^4+24x0^4x3^2x1^2+12x0^6x2^2\\+
%12x3^6x0^2+24x1^2x3^4x2^2+12x0^2x2^6+24x0^2x2^4x3^2\\+
%24x0^4x3^2x2^2+30x0^4x3^4+12x0^6x3^2+3x0^8+30x0^4x2^4\\+
%24x0^2x2^2x3^4+24x1^2x2^4x3^2+3x2^8+12x0^6x1^2+30x2^4x3^4\\+12x0^2x1^6+30x0^4x1^4+12x3^6x2^2+144x1^2x2^2x3^2x0^2+3x1^8\\+12x1^6x2^2+30x1^4x3^4+30x2^4x1^4+12x1^2x3^6+3x3^8+12x2^6x3^2\\+12x2^6x1^2+12x1^6x3^2;\\
%\displaystyle
F_{12}&=&\displaystyle \frac{123}{8}\sum x_i^{12}+\frac{231}{4}\sum x_i^{10}x_j^2+\frac{21}{8}\sum x_i^8x_j^4-\sum \frac{255}{2}\sum x_i^6x_j^6+\frac{949}{2}\sum x_i^8x_j^2x_k^2\\
&&\displaystyle +\frac{1839}{2}\sum x_i^6x_j^4x_k^2+\frac{6111}{4}\sum x_i^4x_j^4x_k^4+1809\sum x_i^6x_j^2x_k^2x_h^2+\frac{7281}{2}\sum x_i^4x_j^4x_k^2x_h^2.
\end{array}
\end{eqnarray*}
Here the sums run over all the indices $i,j,k,h=0,1,2,3$, always being different when appearing together. By applying the map $\phi$, a computer computation with MAPLE shows that
\begin{eqnarray*}
\begin{array}{lll}
\phi(F_6)&=&\displaystyle -\frac{13}{16}\, t_1\cdot t_2,\\
&&\\
\phi(F_8)&=&\displaystyle \frac{3}{64}\, W_1\cdot W_2,\\
&&\\
\phi(F_{12})&=&\displaystyle \frac{3}{256}\, \chi_1\cdot \chi_2
\end{array}
\end{eqnarray*}
as claimed in Corollary \ref{phi}.\\
{\it Proof of Proposition \ref{C1inv}, 2}. The polynomials $F_6,~F_8,~F_{12}$ remain invariant by interchanging $x_2$ with $x_3$, which is what the matrix $C'$ does.
\eprf
{\it Proof of Proposition \ref{quad}}. We write the computations just in the case of the $[3,4,3]$-invariant polynomials. Consider the points $p_1=(i\sqrt{2}:1:1:0)$ and $p_2=(1:i:0:0)$, then $q(p_1)=q(p_2)=0$ and by a computer computation with MAPLE we get $F_6(p_1)=26$, $F_{8}(p_2)=12$ and $F_{12}(p_2)=32$. This shows that $q$ does not divide the polynomials. Since $F_6(p_2)=0$, $F_6$ does not divide $F_{12}$.\eprf 
\begin{rem}
{\rm Observe that an equation for a $[3,4,3]$-invariant polynomial of degree six and for a $[3,3,5]$-invariant polynomial of degree twelve was given by the author in \cite{sarti} by a direct computer computation with MAPLE.}

\end{rem}

%%%%% but computing  the orbit of $p_1$ and $p_2$ under $\tu$ and doing the sums
%\begin{eqnarray*}
%\begin{array}{lll}
%\sum_{g\in\tu}g(f_{11}^{(6)}\cdot f_{22}^{(6)}\cdot\ldots\cdot f_{66}^{(6)})(p_1)&=&\sum_{g\in\tu}(f_{11}^{(6)}\cdot f_{22}^{(6)}\cdot\ldots\cdot f_{66}^{(6)})(g\cdot p_1)\\
%&=&-48(1+i)\\
%\sum_{g\in\tu}g(f_{11}^{(8)}\cdot f_{22}^{(8)}\cdot\ldots\cdot f_{88}^{(8)})(p_2)&=&\sum_{g\in\tu}(f_{11}^{(8)}\cdot f_{22}^{(8)}\cdot\ldots\cdot f_{88}^{(8)})(g\cdot p_2)\\
% &=&12\\
%\sum_{g\in\tu}g(f_{11}^{(12)}\cdot f_{22}^{(12)}\cdot\ldots\cdot f_{1212}^{(12)})(p_2)&=& \sum_{g\in\tu}(f_{11}^{(12)}\cdot f_{22}^{(12)}\cdot\ldots\cdot f_{1212}^{(12)})(g\cdot p_2)\\
%&=&32\\
%\end{array}
%\end{eqnarray*}
\addcontentsline{toc}{section}{  \hspace{0.5ex} References}

%\begin{eqnarray*}
%\begin{array}{l}
%(-1:i\sqrt{2}:0:1), (-1:0:i\sqrt{2}:-1), (0:-1:1:i\sqrt{2})\\
%(i\sqrt{2}:i\sqrt{2}:-i\sqrt{2}+2:i\sqrt{2}+2), (-i\sqrt{2}:-i\sqrt{2}-2:-i\sqrt{2}:i\sqrt{2}+2), (-i\sqrt{2}-2:i\sqrt{2}-2:i\sqrt{2}:i\sqrt{2}), (i\sqrt{2}:-i\sqrt{2}:i\sqrt{2}+2:i\sqrt{2}-2)\\
%(-i\sqrt{2}+2:-i\sqrt{2}-2:-i\sqrt{2}:i\sqrt{2}), (-i\sqrt{2}:-i\sqrt{2}-2:i\sqrt{2}:i\sqrt{2}-2), (i\sqrt{2}+2:-i\sqrt{2}:-i\sqrt{2}+2:i\sqrt{2}), (-i\sqrt{2}+2:-i\sqrt{2}:-i\sqrt{2}-2:-i\sqrt{2})
%\end{array}
%\end{eqnarray*}

%%% Local Variables: 
%%% mode: latex
%%% TeX-master: "artokt02"
%%% End: 


\begin{thebibliography}{EK}
%\bibitem[Be]{benson}
%Benson, D. J.: {\it Polynomial Invariants of Finite Groups}, London Math. Society LNS 190, Cambridge University Press (1993).
\bibitem[Bu]{burnside}
Burnside, W.:{\it Theory of groups of finite order}, Dover Publications, Inc. (1955).
\bibitem[Co1]{coxeter1}
Coxeter, H. S. M.:{\it The product of the generators of a finite group generated by reflections}, Duke Math. J. Vol. 18 (1951) 765-782.
\bibitem[Co2]{coxeter2}
Coxeter, H. S. M.: {\it Regular polytopes (second edition)}, The Macmillan company, New York (1963).
\bibitem[K]{klein}
Klein, F.: {\it Vorlesungen \"uber das Ikosaeder und die Aufl\"osung der Glei\-chungen vom f\"unften Grade}, Nachdr. der Ausg. Leipzig, Teubner 1884, hrsg. mit einer Einf\"uhrung und mit Kommentaren von Peter Slodowy, Birkh\"auser-B. G. Teubner (1993).
\bibitem[Ra]{racah}
Racah, G.: {\it Sulla caratterizzazione delle rappresentazioni irriducibili dei gruppi semisemplici di Lie}, Rend. Acad. Naz. dei Lincei, Classe di Scienze fisiche, matematiche e naturali (8), vol. 8 (1950) 108-112.
\bibitem[Sa]{sarti}
Sarti, A.: {\it Pencils of symmetric surfaces in $\mathbb P^3(\mathbb C)$}, J. of Alg. 246, 429--452 (2001).
\end{thebibliography}
\end{document}